\documentclass[a4paper,reqno,final]{amsart}
\usepackage{latexsym}
\usepackage{amsfonts}
\usepackage{amssymb}
\usepackage{amsmath}
\usepackage{amsthm}
\usepackage{array}
\usepackage{amstext}

\def\R{\mathbb{R}}
\def\C{\mathbb{C}}

\def\ad{\mathrm{ad}}
\def\half{\frac{1}{2}}

\newtheorem{thm}{Theorem}
\newtheorem{defn}[thm]{Definition}
\newtheorem{prop}[thm]{Proposition}
\newtheorem{cor}[thm]{Corollary}

\allowdisplaybreaks[1]

\begin{document}
\title{Penrose limits of homogeneous spaces}
\author{Simon Philip}
\address{School of Mathematics, University of Edinburgh, United
  Kingdom}
\email{S.A.R.Philip@sms.ed.ac.uk} \maketitle

\begin{abstract} We prove that the Penrose limit of a spacetime
along a homogeneous geodesic is a homogeneous plane wave spacetime
and that the Penrose limit of a reductive homogeneous spacetime
along a homogeneous geodesic is a Cahen--Wallach space. We then
consider several homogenous examples to show that these results
are indeed sharp and conclude with a remark about the existence of
null homogeneous geodesics.
\end{abstract}

\section{Introduction}

In \cite{MR56:4680} Penrose introduced a method for taking a
continuous limit of any spacetime to a plane wave. The method
effectively involves ``zooming in" on a null geodesic in such a
way that the metric stays nondegenerate. In \cite{Gueven:2000ru}
G\"uven extended the method to that of supergravity theories where
it is a useful tool for generating new solutions to the
supergravity equations from known ones. Since then several papers
have investigated the properties of these Penrose limits,
\cite{Blau:2003dz}, \cite{Blau:2002dy}, \cite{Blau:2002mw},
\cite{Blau:2003ia}, \cite{Patricot:2003dh}.

Penrose limits have been used as evidence for the $AdS$\textrm{/}
$CFT$ correspondence. The Penrose limits of the $AdS_5 \times S^5$
type $II B$ superstring background were calculated in
\cite{Blau:2002mw}, one of which was shown to be the BFHP
maximally supersymmetric plane wave background \cite{Blau:2002dy}.
String theory in this background is exactly solvable
\cite{Metsaev:2001bj}, \cite{Metsaev:2002re} giving rise to an
explicit form of the $AdS$\textrm{/} $CFT$ correspondence
\cite{Berenstein:2003gb} in which both the gauge theory and the
gravity sides are weakly coupled, allowing many perturbative
checks albeit for a restricted class of observables.

A more general class of background metrics on which string theory
is exactly solvable are the homogeneous plane waves
\cite{Blau:2003rt}, \cite{Papadopoulos:2002bg}. Penrose limits
onto homogeneous plane waves have been investigated, such as the
Penrose limits of the G\"odel-like spacetimes \cite{Blau:2003ia}.
In \cite{Blau:2002mw} it was shown that the dimension of the
isometry algebra never decreases under a Penrose limit. Hence it
seemed a ``natural" assumption that the Penrose limit of a
homogeneous spacetime is always a homogeneous plane wave. However,
in \cite{Patricot:2003dh} it was shown that the cross product of
the homogeneous Kaigorodov spacetime with a sphere has a Penrose
limit which is not itself homogeneous. Consequently the aim of
this paper is give necessary and sufficient conditions on a
spacetime and the null geodesic that guarantee that the Penrose
limit is homogeneous.

Section $2$ gives the definition of the Penrose limit and proves
that this definition is well--defined. We also give a proof of the
covariance property of the Penrose limit stated in
\cite{Blau:2002mw}.

Section $3$ gives some examples of known hereditary properties of
the Penrose limit. Section $4$ contains the background on
homogeneous spaces needed for our results. This includes
descriptions of reductive homogeneous spaces, naturally reductive
homogeneous spaces, the Killing transport and homogeneous
geodesics.

In section $5$ we use the Killing transport to prove that the
Penrose limit of a lorentzian spacetime along a homogeneous
geodesic is a homogeneous plane wave. We then use a similar
approach to prove that the Penrose limit of a reductive
homogeneous spacetime along a homogeneous geodesic is a reductive
homogeneous plane wave.

In section $6$ we use the classification of homogeneous plane
waves that was given in \cite{Blau:2002js} to prove that the
Penrose limit of a reductive homogeneous spacetime along an
absolutely homogeneous geodesic is a naturally reductive plane
wave.

In section $7$ we first show that the Penrose limit of a
non--homogeneous spacetime can be homogeneous. Then we describe
the Kaigorodov spacetime and its Penrose limits as calculated in
\cite{Patricot:2003dh}. We give an example of a Penrose limit of a
reductive homogeneous space along a homogeneous geodesic for which
the homogeneous structure "blows up" but the limiting spacetime is
still homogeneous. We also give an example of a Penrose limit of a
non--reductive homogeneous spacetime along a homogeneous geodesic
which is still non--reductive homogeneous.

Finally in section $8$ we show that while there must exist at
least one homogeneous geodesic in any reductive homogeneous
spacetime \cite{MR2002b:53076}, \cite{MR2001f:53104}, there may
not exist any null absolutely homogeneous geodesics.

\section{What is a Penrose limit?}

Let $(M,g)$ be a smooth ($n+1$)-manifold with a lorentzian metric.
Let $\gamma$ be a null geodesic of $(M,g)$. Then given a point $x
\in \gamma$ there exists a coordinate neighborhood $(U, \mu)$,
$\mu : U \to \R^{n+1}$, of $x$ defining coordinates $\mu(y) =
(u(y),v(y),[y^k(y)])$, where $u$ is a coordinate along $\gamma$,
such that in $U$ the metric may be written as

\begin{equation} \label{2} g = du dv + \alpha dv^2 + \sum_{i=1}^{n-1}
\beta_i dy^i dv + \sum_{i,j=1}^{n-1} C_{ij} dy^i dy^j.
\end{equation} Here $\alpha, \beta_i, C_{ij}$ are functions of
$(u,v,[y^k])$ and $(C_{ij})$ is positive definite.

To choose such coordinates one chooses a one-parameter family
 of hypersurfaces parameterized by $v$ and foliated by null geodesics.
The coordinate along the prescribed geodesics is given by $u$ and
$\gamma$ is given by $(u,0,0)$.

In other words, one chooses a local extension of the null geodetic
tangent vectorfield $\frac{\partial}{\partial u}$ of $\gamma$ to a
null geodetic vector field in a neighborhood of $x$. Then one
chooses $(n-1)$-submanifolds on which the restricted metric is
Riemannian and allows $v$ to be the parameter labelling these
submanifolds.

Let $\Omega \in (0, \infty)$. Consider the map

\begin{equation} \begin{split} \psi_{\Omega} & : \R^{n+1} \to
\R^{n+1}\\ & : (u,v,[y^k]) \mapsto (u, \Omega^2 v,  [ \Omega
y^k]).
\end{split}
\end{equation} This map induces the following change of coordinates:

\begin{equation} \phi^U_{\Omega}= \mu^{-1} \circ \psi_{\Omega} \circ \mu :
U \to U . \end{equation} (If necessary, to make this well defined,
we may need to shrink $U$ so that it does not contain any
``holes".) By patching together such coordinate neighborhoods
along $\gamma$ we may think of $\phi_{\Omega}$ as a diffeomorphism
from a tubular neighborhood of $\gamma$ to a tubular
subneighborhood. Then we define the \textbf{Penrose limit} of
$(M,g,\gamma)$ to be this tubular neighborhood together with the
metric

\begin{equation} \label{7}
  \begin{split}
    g_{Pl} & = \lim_{\Omega \to 0} \Omega^{-2} (\phi_{\Omega}^{-1})^*
    g\\
    & = du dv + \sum_{i,j=1}^{n-1} C_{ij}(u,0,0) dy^i dy^j~.
  \end{split}
\end{equation}
Notice that at $\Omega = 0$, $\phi_{\Omega}$ is no longer a
diffeomorphism.

\begin{prop}
  $g_{Pl}$ is defined independently of choice of coordinates
  putting $g$ in the form \eqref{2}.
\end{prop}

\begin{proof}
  Let $(r,s,[x^i])$ be a different choice of coordinates such that
  \begin{equation} \label{6}
    g = dr ds + \rho ds^2 + \sum_{i=1}^{n-1} \psi_i dx^i ds +
    \sum_{i,j=1}^{n-1} \Theta_{ij} dx^i dx^j~,
  \end{equation}
  where $\rho, \psi_i, \Theta_{ij}$ are functions of $(r,s,[x^i])$ and
  $(\Theta_{ij})$ is positive definite. As both $u$ and $r$ are
  parameters along the geodesic $\gamma$ we may as well choose them
  equal $u=r$. An easy check shows that the change of coordinates
  matrix must be of the form
  \begin{equation*}
    \begin{pmatrix}
      dr \\ ds \\ dx^i
    \end{pmatrix} =
    \begin{pmatrix}
      1 & 0 & 0 \\ 0 & 1 & 0 \\ 0 & c^i & e^i_k
    \end{pmatrix}
    \begin{pmatrix}
      du \\ dv \\ dy^k
    \end{pmatrix}
  \end{equation*}
  and that under this
  \begin{equation}
    \Theta_{ij} e^i_k e^j_l = C_{kl}.
  \end{equation}
  In fact $c^i$ must also be zero because the second row in the matrix
  equation above shows that $s = v+K$, $K$ a constant, and the change of
  basis matrix for the dual basis to the one-forms above is the
  inverse transpose:
  \begin{equation*}
    \begin{pmatrix}
      \partial u \\ \partial v \\ \partial y^i
    \end{pmatrix}
    =
    \begin{pmatrix}
      1 & 0 & 0 \\
      0 & 1 & -c^i (e^{-1})^i_k \\
      0 & 0 & (e^{-1})^i_k
    \end{pmatrix}
    \begin{pmatrix}
      \partial r \\ \partial s \\ \partial x^k
    \end{pmatrix}
  \end{equation*}
  As $e^i_k$ is nondegenerate we must have $c^i=0$. Putting this
  into the Penrose limit metric \eqref{7}
  \begin{equation}
    \begin{split}
      drds + \sum_{i,j=1}^{n-1} \Theta_{ij}(r,0,0) dx^i dx^j & = dudv
      +
      \sum_{i,j,k,l=1}^{n-1}\Theta_{ij}(r,0,0) e^i_k e^j_l  dy^k dy^l \\
      &= dudv + \sum_{k,l=1}^{n-1} C_{kl}(u,0,0) dy^k dy^l~.
  \end{split}
\end{equation}
\end{proof}

In the recent paper \cite{Blau:2003dz} a covariant description of
the Penrose limit without reference to the adapted coordinates is
given.

A sufficient condition for telling when two Penrose limits will be
isometric is the following (the statement of this Theorem appeared
in \cite{Blau:2002mw} although the proof did not).

\begin{thm}[Covariance of the Penrose limit]
  Let $(M,g),(M',g')$ be two lorentzian manifolds. Let $\gamma$ and $
  \gamma' $ be two null geodesics inside $M$ and $M'$ respectively.
  Let $f: M_{\gamma} \to M_{\gamma'}'$ be an isometry of tubular
  neighborhoods of $\gamma$ and $ \gamma' $ which maps $\gamma$ onto
  $\gamma' $. Then the Penrose limits of $(M,g)$ and $(M',g')$ along
  $\gamma$ and $\gamma'$ respectively are isometric.
\end{thm}

\begin{proof}
  Let $(U, \mu=(u,v,[y^k]))$ be a coordinate neighborhood of a point
  $x$ on $\gamma$ such that the metric $g$ takes the form \eqref{2}.
  Define a coordinate neighborhood $(f(U), \mu'=(u',v',[{y'}^k]))$
  about $f(x)$ by
  \begin{equation}
    \mu'(f(x)) = \mu (x),
  \end{equation}
  so that $u' = u \circ f^{-1}$ is a coordinate along $\gamma'$. As $g
  = f^*g'$, then $g'$ also takes the form of \eqref{2} in this
  neighborhood.

  Now consider $f \circ \phi_{\Omega}^U: U \to U'$. We have
  \begin{equation}
    \begin{split}
      f \circ \phi_{\Omega}^U &= f \circ \mu^{-1} \circ \psi_{\Omega}
      \circ \mu \\
      &= f \circ (\mu' \circ f)^{-1} \circ \psi_{\Omega} \circ (\mu'
      \circ f) \\
      &= {\mu'}^{-1} \circ \psi_{\Omega} \circ \mu' \circ f\\
      &= \phi_{\Omega}^{U'} \circ f~.
    \end{split}
\end{equation}
Therefore
\begin{equation}
  \begin{split}
    g_{Pl} &= \lim_{\Omega \to 0} \Omega^{-2} (\phi_{\Omega}^U)^* g
    \\
    & = \lim_{\Omega \to 0} \Omega^{-2}(\phi_{\Omega}^U)^*
    f^* g' \\
    &= \lim_{\Omega \to 0} \Omega^{-2}(f \circ \phi_{\Omega}^U)^* g'
    \\
    &= \lim_{\Omega \to 0} \Omega^{-2} (\phi_{\Omega}^{U'}
    \circ f)^* g'\\
    &= \lim_{\Omega \to 0} \Omega^{-2} f^* \circ
    (\phi_{\Omega}^{U'})^* g'\\
    & = f^* g'_{Pl}~.
  \end{split}
\end{equation}

\end{proof}

\section{Hereditary properties}

We say that a property of the metric $g$ is \textbf{hereditary} if
the Penrose limit $g_{Pl}$ has the same property. For example,

\begin{prop} Suppose $(M,g)$ is locally symmetric/conformally flat. Then
 $(M_{\gamma}, g_{Pl})$ is locally symmetric/conformally flat. If
$(M,g)$ is Einstein then $(M_{\gamma}, g_{Pl})$ is Ricci flat, in
particular it is Einstein.
\end{prop}
\begin{proof}
  Let $\nabla_{\Omega}, R_{\Omega}$ denote the connection and
  curvature of $g_{\Omega} := \Omega^{-2} (\phi^{-1}_{\Omega})^* g$
  respectively. As $\phi_{\Omega}$ is a diffeomorphism if $\nabla R =
  0 $ then $\nabla_{\Omega} R_{\Omega} = 0$ for $\Omega > 0$. By a
  continuity argument we see that $\nabla_{Pl} R_{Pl} = 0$.\\
  If $Ric(g) = \lambda g$ then
  \begin{equation}
    Ric(g_{\Omega}) = Ric(\Omega^{-2} (\phi^{-1}_{\Omega})^* g ) =
    Ric((\phi^{-1}_{\Omega})^* g ) = \lambda (\phi^{-1}_{\Omega})^*
    g~.
  \end{equation}
  This gives
  \begin{equation}
    Ric(g_{\Omega}) = \Omega^2 \lambda g_{\Omega}~,
  \end{equation}
  and by continuity we see that $Ric(g_{Pl}) = 0$.
\end{proof}

These hereditary properties can be used to easily compute the
Penrose limits of anti de~Sitter space $AdS$. Anti de~Sitter space
is Einstein and conformally flat hence any Penrose limit is Ricci
flat and conformally flat and thus flat.

In \cite{Blau:2002mw} the case of $AdS \times S$ is considered. It
is a symmetric space and is shown to have two non-isometric null
geodesics leading to two non-isometric Penrose limits which are
flat space and a symmetric plane wave.

Another useful hereditary property is that of geodesic
completeness:

\begin{thm} Suppose $(M,g)$ is a geodesically complete lorentzian
manifold. Then the Penrose limit along any null geodesic is
geodesically complete.
\end{thm}

\begin{proof} Let $\gamma(t)$ be a geodesic with respect to $\nabla_{Pl}$ for $t \in
[a,b]$. Without loss we may assume that $\gamma$ is contained in a
normal coordinate neighborhood of some point on $\gamma$ so that
there is a unique geodesic from $\gamma(a)$ to $\gamma(b)$ with
respect to $\nabla_{\Omega}$ for $\Omega \in [0,1]$ (which is
possible because $\nabla_{\Omega}$ varies continuously with
respect to $\Omega$ and $[0,1]$ is compact.) . Let
$\gamma_{\Omega}$ be the unique geodesic with respect to
$\nabla_{\Omega}$ between $\gamma(a)$ and $\gamma(b)$. Then
$\gamma_{\Omega}(t) $ may be extended to $(-\infty, \infty)$ as
$\nabla_1$ is geodesically complete and $\phi_{\Omega}$ is a
diffeomorphism. Continuity implies that the sequence of geodesics
$\gamma(\Omega)$ for $\Omega= \frac{1}{k}$ ``converges" to
$\gamma$ in the following sense. Any neighborhood of any point on
$\gamma$ intersects all but a finite number of geodesics of the
sequence. Therefore, by continuity of the geodesic equation with
respect to $\Omega$, we have that $\gamma$ may be extended beyond
$(a,b)$.
\end{proof}

One last hereditary property, as noted in the introduction, is the
following:

\begin{prop} The dimension of the isometry algebra of $g_{Pl}$ is no less
than the dimension of the isometry algebra of $g$.
\end{prop}
\begin{proof} See \cite{MR40:3875} or \cite{Blau:2002mw}. \end{proof}

\section{Homogeneous spaces and homogeneous structures}

In this section we will give the definitions and results we need
in relation to homogeneous spaces.

\begin{defn} A connected lorentzian manifold $(M,g)$ is
  \textbf{homogeneous} if its group of isometries acts transitively on
  $M$.
\end{defn}

When this is the case then $M$ can be written $M= G/H$ where $G$
is the group of isometries and $H$ is a closed subgroup.

\begin{defn} A homogeneous space $M= G/H$ is \textbf{reductive} when there
exists a subspace $\mathfrak{m} \cong T_pM \subset\mathfrak{g}$
such that

\begin{enumerate}
\item $\mathfrak{g} = \mathfrak{h} \oplus \mathfrak{m},$
\item $[\mathfrak{h}, \mathfrak{m}] \subset \mathfrak{m}.$
\end{enumerate} It is \textbf{symmetric} if it also satisfies
\begin{equation} [\mathfrak{m}, \mathfrak{m}] \subset
\mathfrak{h}. \end{equation}
\end{defn} (In fact, strictly speaking this is the definition
of \textbf{weakly reductive}. However for the rest of this paper
we shall assume that $H$ is connected, in which case they are the
same thing.)

\begin{defn} Let $o$ denote the coset of $H$ in $M$ and fix a frame
$u_0:\R^n \to T_oM$ of the frame bundle $F$. Define the
\textbf{linear isotropy representation} $\lambda:H \to GL(n,\R)$
by

\begin{equation} \lambda (h):= u_o^{-1} \circ h_* \circ u_o,
\end{equation} where $h \in H$, $h_*: T_oM \to T_oM$ denotes the
differential of $h$ at $o$.

\end{defn}

\begin{thm} Let $F$ be the frame bundle of $M=G/H$ a reductive
homogeneous space of dimension $n$ with decomposition $
\mathfrak{g} = \mathfrak{h} \oplus \mathfrak{m}$. Then there is a
one-to-one correspondence between the set of $G$-invariant
connections in $F$ and the set of linear maps
$\Lambda_{\mathfrak{m}} : \mathfrak{m} \to \mathfrak{gl}(n, \R)$
such that

\begin{equation} \Lambda_{\mathfrak{m}}(\ad h(X)) = \ad( \lambda(h))
(\Lambda_{\mathfrak{m}}(X)), \end{equation} for $X \in
\mathfrak{m}$ and $h \in H$.

The correspondence is given by

\begin{equation} \omega_{u_o}(\tilde{X}) = \begin{cases} \lambda(X)
& \textrm{if} \; X \in \mathfrak{h}, \\
\Lambda_{\mathfrak{m}}(X) & \textrm{if} \; X \in \mathfrak{m},
\end{cases}
\end{equation} where $\omega$ is the connection one-form,
$\tilde{X}$ is the natural lift of $X \in \mathfrak{g}$ to $F$ and
$\lambda$ is not only as above $H \to GL(n,\R)$ but also the
induced Lie algebra homomorphism $\mathfrak{h} \to
\mathfrak{gl}(n,\R)$.

\end{thm}
\begin{proof} See chapter X, Theorem $2.1$ in \cite{MR38:6501}. \end{proof}

\begin{defn} The connection obtained by taking $\Lambda_{\mathfrak{m}}
= 0 $ is called the \textbf{canonical connection}. \end{defn}

The canonical connection can also be described in the following
way. Let $\theta$ be the left-invariant Maurer--Cartan form of $G$

\begin{equation} \theta_g(X) := (L_g)^*(X), \end{equation} where
$L_g$ denotes left multiplication by $g$ and $*$ denotes
differentiation. Let $\sigma: U \to G$ be a local coset
representative. Then the pull back of $\theta$ by $\sigma$ splits
as

\begin{equation}
  \sigma^*(\theta) = \theta^{\mathfrak{h}} + \theta^{\mathfrak{m}}~,
\end{equation}
where $\theta^{\mathfrak{h}}_x(X) \in \mathfrak{h},
\theta^{\mathfrak{m}}_x(X) \in \mathfrak{m}$.  The one-form
$\theta^{\mathfrak{h}}$ defines the connection one-form for the
canonical connection.

The geodesics of the canonical connection are curves $\gamma(t)$
of the form
\begin{equation}
  \exp( tX), \; t \in \R, X \in
  \mathfrak{g}~.
\end{equation}
If $(M,g)$ is symmetric then the canonical connection coincides
with the Levi--Civita connection.

\begin{thm}
  \label{canconn}
  The canonical connection of a reductive homogeneous space is
  complete.
\end{thm}

\begin{proof} See chapter X, Corollary $2.5$ in \cite{MR38:6501}. \end{proof}

\begin{thm}[(\cite{MR21:1628},\cite{MR22:1863},\cite{MR93h:53047}]
  \label{3}
  Let $(M,g)$ be a reductive lorentzian homogeneous manifold with
  Levi--Civita connection $\nabla$. Then there exists a $(2,1)$ tensor
  $T$ defining a metric connection $\tilde{\nabla} := \nabla - T$ with
  curvature $R$ such that $\tilde{\nabla} T = \tilde{\nabla} R = 0$.
\end{thm}

\begin{proof} Write $M = G/H$, with decomposition $ \mathfrak{g} =
\mathfrak{h} \oplus \mathfrak{m}$. Let $\tilde{\nabla}$ be the
canonical connection of $M$. Let $T = \nabla - \tilde{\nabla}$. As
$G$ acts by isometries, $\nabla$ is also $G$-invariant. Hence $T$
and $R$ are $G$-invariant. Therefore, see \cite{MR38:6501}, they
are both parallel with respect to $\tilde{\nabla}$.
\end{proof}

(The first version of this Theorem for riemannian signature
appeared in \cite{MR21:1628}. This was re-interpreted in terms of
the canonical connection in \cite{MR22:1863} and extended to the
pseudo-riemannian case in \cite{MR93h:53047}.)

\textbf{Remarks}: \begin{enumerate} \item $T$ is not necessarily
the torsion of $\tilde{\nabla}$ (and not necessarily
skew-symmetric in it lower indices.) If $\tilde{\Gamma}^i_{jk}$
are the Christofel symbols of $\tilde{\nabla}$ and $\Gamma^i_{jk}$
of $\nabla$ and $\tau$ is the torsion of $\tilde{\nabla}$ then
\begin{equation}
  \begin{split}
    \tau^i_{jk} &= \tilde{\Gamma}^i_{jk} - \tilde{\Gamma}^i_{kj}\\
    &= \tilde{\Gamma}^i_{jk} -\Gamma^i_{jk} + \Gamma^i_{kj} -
    \tilde{\Gamma}^i_{kj} \\
    &= -T^i_{jk} +T^i_{kj}~.
  \end{split}
\end{equation} i.e. $\tau$ is the skew-symmetrization of $T$. In fact
$\tau(X,Y)|_{\mathfrak{m}} = -[X,Y]_{\mathfrak{m}}$, the component
of $[X,Y]$ lying in $\mathfrak{m}$ where $X,Y \in \mathfrak{m}$
(see chapter X, Theorem $2.6$ in \cite{MR38:6501}.) Also the
restriction of $T$ to $\mathfrak{m}$ is given by
\begin{equation}
  T(X,Y)|_{\mathfrak{m}} = \half
  [X,Y]_{\mathfrak{m}} + U(X,Y),
\end{equation}
where $U$ is the symmetric bilinear mapping of $\mathfrak{m}
\times \mathfrak{m}$ into $\mathfrak{m}$ defined by
\begin{equation}
  2g_o(U(X,Y),Z) =
  g_o(X,[Z,Y]_{\mathfrak{m}}) +
  g_o([Z,X]_{\mathfrak{m}},Y),
\end{equation}
where $X,Y,Z \in \mathfrak{m}$ (see chapter X, Theorem $3.3$ in
\cite{MR38:6501}.)

\item If $(M,g)$ is symmetric then $\nabla R = 0$ and we can take $T=0$.
\end{enumerate}
Such a $T$-tensor is called a \textbf{homogeneous structure}. A
given homogeneous manifold $M$ may have many different homogeneous
structures. Each corresponding to a different choice of groups $G$
and $H$. For example, the $7$-sphere $S^7 = SO(8)/SO(7) =
Spin(7)/G_2=Sp(2)/Sp(1)$. (For a review of Penrose limits from the
point of view of homogeneous structures see
\cite{FigueroaMeessenPhilip}.)

\begin{defn}
  $(M,g)$ is called \textbf{naturally reductive} if there exists a
  homogeneous structure $T$ with $U =0$, i.e. if $\tau =
  T$.
\end{defn}
While reductivity is a property of the isotropy representation,
natural reductivity is also a property of the metric.

\begin{prop} \label{5}
  Let $(M,g)$ be naturally reductive. Then the geodesics of the
  Levi--Civita connection coincide with the geodesics of the canonical
  connection.
\end{prop}

\begin{proof}
  Let $\nabla$ denote the Levi--Civita connection and $\tilde{\nabla}$
  the canonical connection corresponding to the homogeneous structure
  $T$ with $U=0$. Then $T$ will be skew-symmetric in its lower
  indices and consequently
  \begin{equation}
    \nabla_{X} X
    =\tilde{\nabla}_{X}X +
    T(X,X) =
    \tilde{\nabla}_{X}X~.
  \end{equation}
  Hence the geodesic equations for the two connections are the same.
\end{proof}

Theorem \ref{3} can be rewritten and then a converse constructed:

\begin{thm}
  \label{1}
  Let $(M,g)$ be a connected, simply connected, lorentzian manifold.
  Then $(M,g)$ is reductive homogeneous if and only if there exists a
  complete affine metric connection $\tilde{\nabla}$ with torsion
  $\tau$ and curvature $R$ such that $\tilde{\nabla} \tau =
  \tilde{\nabla} R = 0$.
\end{thm}

\begin{proof} See chapter X, Theorems $2.6$--$2.8$ in \cite{MR38:6501}. \end{proof}

Above we have Theorems \ref{3} and \ref{1}, which describe the
reductive homogeneous property in terms of a metric connection on
the tangent bundle. In fact we can describe the Killing vectors of
an arbitrary pseudo-riemannian manifold as parallel vector fields
of a covariant derivative on an extended bundle.

Let $X$ be a vector field on a lorentzian manifold $(M,g)$. Let
$A_{X} : Y \mapsto -\nabla_{Y}X$. Then $X$ is a Killing vector if
and only if $A_{X}$ is skew-symmetric with respect to $g$. As a
consequence of the Killing identity we also have the equation:
\begin{equation*}
  \nabla_{X} A_{\zeta} = -R(X, \zeta)
\end{equation*}

Now consider the bundle $\mathfrak{E} = TM \oplus
\mathfrak{so}(TM)$. If we define a covariant derivative $D$ on
$\mathfrak{E}$ by
\begin{equation*}
  D_{X} (\zeta, A) := (\nabla_{X} \zeta + A(X), \nabla_{X}A +
  R(X, \zeta)).
\end{equation*} Then the parallel sections of $\mathfrak{E}$ with respect to $D$
are precisely the Killing vectors of $g$. Thus a Killing vector is
completely determined by
\begin{equation*}
  (\zeta(p), A_{\zeta}(p))
\end{equation*}
at any point $p$ and by parallel translation by the covariant
derivative $D$.

Finally we make the

\begin{defn}
  A geodesic $\gamma$ is called homogeneous if it is the orbit of a
  $1$-parameter subgroup of isometries.
\end{defn}

\textbf{Note}: on a riemannian space this definition is equivalent
to writing the geodesic in the form $\gamma( t) = exp(tX)_o$ for
some $X \in \mathfrak{g}$ (see \cite{MR92m:53084}.) However, if
the geodesic $\gamma( t)$ is null one may have to change its
parameterization in order to write it in the form $exp(sX)_o$. We
call a geodesic of the form $\gamma( t) = exp(tX)_o$ an
\textbf{absolutely homogeneous geodesic}. Also notice that a
homogeneous geodesic is not necessarily a geodesic of the
canonical connection as a geodesic of the canonical connection is
of the form $\exp (tX)$ with $X \in \mathfrak{m}$. We shall call a
homogeneous geodesic which is a geodesic of the canonical
connection a \textbf{canonical homogeneous geodesic}.

These will be the geodesics of interest when deciding whether a
Penrose limit is homogeneous or not. A useful criteria for
distinguishing homogeneous geodesics is the following,

\begin{prop} \label{geo}
  Suppose $M$ is a lorentzian reductive homogeneous space. The
  geodesic $\gamma(t)$ with $\gamma(0) = o$ and
  $ \gamma'(0) = X \in \mathfrak{g}$ is a homogeneous geodesic if and only
  if \begin{equation} \label{geovec} B(X_{\mathfrak{m}}, [Z,X]_{\mathfrak{m}}) =
  \lambda B(X_{\mathfrak{m}},Z_{\mathfrak{m}}) \end{equation} for
  all $Z \in \mathfrak{g}$ and some $\lambda \in \R$. It is
  absolutely homogeneous if and only if $\lambda = 0$.
\end{prop}
\begin{proof} This is a slight generalization of the proof given
in \cite{MR92m:53084}. \end{proof}

\begin{defn} A vector $X \in \mathfrak{g}$ which satisfies
\eqref{geovec} is called a \textbf{geodesic vector}. \end{defn}

\textbf{Note}: by putting $Z=X$ in \eqref{geovec}, we see that if
$B$ has riemannian signature then we must have $\lambda = 0$.

\textbf{Remark}: Suppose $\gamma$ is a geodesic parameterized by
$u$. If $\gamma$ is homogeneous then there is a Killing vector
$\zeta$ such that $ \zeta_p = \gamma'_p$ at all points $p \in
\gamma$. But the geodesic vector field $\frac{\partial}{\partial
u}$ is not necessarily a Killing vector field. If the Killing
vector field $\zeta$ is given by $\frac{\partial}{\partial t^1}$,
then $t^1$ may not be part of a \textbf{twist--free coordinate
system}; that is a coordinate system $(t^1, \dots, t^n)$, in which
we can write the metric in the form $g_{ij} dt^i dt^j$ such that
$d(g_{1i} dt^i) =0$. In particular, if $\gamma$ is a null
homogeneous geodesic then we may not be able to write $g$ in the
form of \eqref{2} with $\frac{\partial}{\partial u}$ a Killing
vector.

\begin{prop} A geodesic $\gamma$ is homogeneous if and only if
there exists a solution $(\gamma, A)$ to the Killing transport
equations with $A(\gamma') = 0$.
\end{prop}

\section{Penrose limits along homogeneous geodesics}

In this section we will give three Theorems which give sufficient
conditions for homogeniety to be hereditary.

\begin{thm} \label{killingpl}
  The Penrose limit along a null geodesic $\frac{\partial}{\partial
    u}$ which is a Killing vector is flat.
 \end{thm}

\begin{proof}
  \begin{equation*}
    \begin{split}
      0 = \mathfrak{L}_{\frac{\partial}{\partial u}} g &=
      d(i_{\frac{\partial}{\partial u}}g) + i_{\frac{\partial}{\partial
          u}} dg \\
      &= d(dv) + \frac{\partial \alpha}{\partial u} dv^2 +
      \frac{\partial \beta_i}{\partial u} dv dy^i + \frac{\partial
        C_{ij}}{\partial u} dy^i dy^j~.
    \end{split}
  \end{equation*}
  Therefore $C$ is independent of $u$ and hence $g_{Pl}$ is flat.
\end{proof}

\begin{thm} \label{plhomog} The Penrose limit of a lorentzian metric along a homogeneous
geodesic $\gamma$ is homogeneous. \end{thm}

\begin{proof}
  On a plane wave $ds^2 = dudv + C_{ij}dy^idy^j$, we have the Killing vectors
  $$\frac{\partial}{\partial v}, \frac{\partial}{\partial y^1}, \dots,
  \frac{\partial}{\partial y^{n-1}}$$  which are independent at each
  point $p$. So to prove $g_{Pl}$ is locally homogeneous it is enough to show
  it has a Killing vector which agrees with $\frac{\partial}{\partial
    u}= \gamma'$ at $p$.\\
  Suppose that $\zeta$ is a Killing vector such that
  $\zeta |_{\gamma} = f \frac{\partial}{\partial u}|_{\gamma}$.
  Then $\zeta|_{\gamma}$ is generated by Killing transport of
  $(\zeta(p), A_{\zeta}(p))$ along $\gamma$. Now by
  definition,
  \begin{equation*}
    \left( A_{\zeta} f \gamma'
    \right)|_{\gamma} = \left( A_{\zeta} \zeta \right)|_{\gamma}  =
    0,
  \end{equation*}
  where by $|_{\gamma} $ we mean restriction to $\gamma \in
  M$,
  not restriction of the tangent bundle. Therefore, if we write $A_{\zeta}$ in components:
  \begin{equation*}
    A_{\zeta} = \sum_{i,j} \left(A_{\zeta}
    \right)^j_{i} dx^i \otimes \frac{\partial}{\partial x^j},
  \end{equation*}
  we see that
  \begin{equation*}
    \left( A_{\zeta} \right)^{y^i}_u |_{\gamma} = \left(
      A_{\zeta} \right)^v_u |_{\gamma} = 0.
  \end{equation*}
  Also, as $\zeta$ is a Killing vector, we have
  \begin{equation*}
    g\left(A_{\zeta} \frac{\partial}{\partial y^i} ,
      \frac{\partial}{\partial u}\right)|_{\gamma} = -g\left(
      \frac{\partial}{\partial y^i} , A_{\zeta}
      \frac{\partial}{\partial u}\right)|_{\gamma} = 0~.
\end{equation*}
Therefore,
\begin{equation*}
  \left( A_{\zeta} \right)^v_{y^i} |_{\gamma} = 0~.
\end{equation*}

Now consider the pull-back of the Killing transport covariant
derivative under the Penrose limit map $\phi_{\Omega}$;
\begin{equation*}
  (\phi_{\Omega}^{-1})^* D_{\zeta}(X)|_{\gamma}
  =(\phi_{\Omega}^{-1})^*\nabla_{\zeta}(X)
  |_{\gamma} - (\phi_{\Omega}^{-1})^*
  A_{\zeta}(X)|_{\gamma}~.
\end{equation*}
The components of $A_{\zeta}|_{\gamma}$ scale under the Penrose
limit map in the following way:
\begin{align*}
  \left( A_{\zeta} \right)^{y^i}_u |_{\gamma} \mapsto &
  \Omega^{-1} \left( A_{\zeta} \right)^{y^i}_u |_{\gamma} \\
  \left( A_{\zeta} \right)^v_u |_{\gamma} \mapsto &
  \Omega^{-2} \left( A_{\zeta} \right)^v_u |_{\gamma} \\
\left( A_{\zeta} \right)^v_{y^i} |_{\gamma} \mapsto & \Omega^{-1}
\left( A_{\zeta} \right)^v_{y^i} |_{\gamma}
\end{align*}
and other components which either stay constant or tend to zero as
$\Omega \to 0$. Taking the limit as $\Omega \to 0$ we see from
above that the three components of $A_{\zeta}$ that could ``blow
up" are in fact zero. Therefore
\begin{equation*}
  (D_{Pl})_{\zeta}(X)(u,v,y) := \lim_{\Omega
    \to 0} [D_{\zeta}(X)(u,0,0)]
\end{equation*}
is well-defined and along with
\begin{equation*}
  (D_{Pl})_{\zeta}(A) :=
  (\nabla_{Pl})_{\zeta} A - R_{Pl}(\zeta,
  X)~,
\end{equation*}
defines a Killing transport covariant derivative on along $\gamma$
with respect to $g_{Pl}$.\\
Therefore, parallel translation by $D_{Pl}$ along $\gamma$
generates the remaining Killing vector needed.
\end{proof}

\begin{cor} \label{plhomogcor} If $\gamma$ is an (absolutely) homogeneous geodesic of $g$ then it
is also an (absolutely) homogeneous geodesic of the Penrose limit
of $g$ along $\gamma$.
\end{cor}

When $(M,g)$ is a reductive homogeneous manifold we can use the
same strategy as above to construct a homogeneous structure on the
Penrose limit:

\begin{prop} \label{red} The Penrose limit of a reductive
lorentzian homogeneous manifold along a canonical homogeneous
geodesic is locally reductive homogeneous.
\end{prop}

\begin{proof}
  Let $(M,g)$ be a reductive homogeneous space with a null homogeneous
  geodesic $\gamma$. From the Ambrose--Singer Theorem we have a
  connection $\tilde{\nabla}$ such that $\tilde{\nabla} T =
  \tilde{\nabla} R = 0$. Let $M_{\gamma}$ be a tubular neighborhood of
  $\gamma$ as above and consider $\phi_{\Omega}( M_{\gamma})$.  Now
  $\phi_{\Omega}$ is a diffeomorphism for $\Omega \neq 0$ so
  $\phi_{\Omega}( M_{\gamma})$ is reductive homogeneous for $\Omega
  >0$. This defines the metric connection
  \begin{equation} \label{4}
    \tilde{\nabla}_{\Omega} := (\phi_{\Omega}^{-1})^* \tilde{\nabla} =
    (\phi_{\Omega}^{-1})^* \nabla - (\phi_{\Omega}^{-1})^* T~.
  \end{equation}
  $\gamma$ is a homogeneous geodesic of $(M,g)$ so $T$ is a tensor of
  type $(2,1)$
  \begin{equation}
    T = \nabla - \tilde{\nabla} = \sum_{i,j,k=1}^{n+1} T_{ij}^k dx^i
    \otimes dx^j \otimes \frac{\partial}{\partial x^k}~.
  \end{equation}
  Under the Penrose limit map $\phi_{\Omega}$ the coefficients scale in
  the following way
  \begin{equation*}
    T_{uy^i}^v \mapsto \Omega^{-1} T_{uy^i}^v |_{\gamma} \qquad
    T_{uu}^v \mapsto \Omega^{-2} T_{uu}^v |_{\gamma}\qquad
    T_{uu}^{y^i} \mapsto \Omega^{-1} T_{uu}^{y^i} |_{\gamma}~,
  \end{equation*}
  and terms which either remain the same or tend to $0$ in the limit
  $\Omega \to 0$.\\
  Suppose that $\gamma$ is a canonical homogeneous geodesic. Then there is a
  Killing vector $\zeta$ such that $\zeta
  |_{\gamma} = f \frac{\partial}{\partial u}|_{\gamma}$. Then
  $\zeta|_{\gamma}$ is generated by parallel transport by
  the canonical connection of $\zeta(p)$ along $\gamma$.
  Now by definition,
  \begin{equation*}
    (\nabla_{\gamma'} \gamma')|_{\gamma} = 0 \quad \textrm{ and } \quad (\tilde{\nabla}_{\zeta} \zeta)|_{\gamma} =0
  \end{equation*}
  where by $|_{\gamma} $ we mean restriction to $\gamma \in M$ not
  restriction of the tangent bundle. Thus \begin{equation*} 0 =(\tilde{\nabla}_{f\gamma} f
  \gamma)|_{\gamma}=(\nabla_{f\gamma'} f\gamma')|_{\gamma} - T(f
  \gamma', f \gamma')|_{\gamma} = f df(\gamma') \gamma'|_{\gamma} -
  f^2T(\gamma',\gamma')|_{\gamma}, \end{equation*} and therefore,
  \begin{equation*}
    T^{y^i}_{uu} |_{\gamma} = T^v_{uu} |_{\gamma} = 0.
  \end{equation*}
  Also, as $\tilde{\nabla}$ is metric
  we have
  \begin{equation}
    \begin{split}
      0 &= (\tilde{\nabla}_{W}
      g)(X,Y)\\
      &= (\nabla_{W} g)(X,Y) +
      g(T_{W}X,Y) +
      g(X,T_{W}Y)\\
      &= g(T_{W}X,Y) +
      g(X,T_{W}Y)~,
    \end{split}
  \end{equation}
  as $\nabla$ is metric. Hence using \eqref{2} we see that
  \begin{equation}
    T_{uy^i}^v|_{\gamma} =  0~.
  \end{equation}
  The Levi--Civita connection of the Penrose limit along $\gamma$,
  $\nabla_{Pl}$, is equal to
  \begin{equation}
    \lim_{\Omega \to 0} (\phi_{\Omega}^{-1})^* \nabla~.
  \end{equation}
  Also above shows that the limit $T_{Pl}|_{\gamma} := \lim_{\Omega
    \to 0}(\phi_{\Omega}^{-1})^* T|_{\gamma} $ is well defined on
  $\gamma$.  Thus, by \eqref{4}, the limit  $
  \tilde{\nabla}_{Pl}|_{\gamma} := lim_{\Omega \to
    0}\tilde{\nabla}_{\Omega}|_{\gamma}$ is well defined. Now
  \begin{equation*}
    \{
    \tilde{\nabla}_{\Omega}|_{\gamma} g_{\Omega}|_{\gamma} \; | \;
    \Omega \in [0,1] \}
  \end{equation*}
  is a continuous path in the space of tensors
  of type $(3,0)$ on $\gamma$. Therefore continuity shows $
  \tilde{\nabla}_{Pl}|_{\gamma} g_{Pl}|_{\gamma} = 0$.  Similarly
  \begin{equation}
    \tilde{\nabla}_{Pl}|_{\gamma} g_{Pl}|_{\gamma} = \tilde{\nabla}_{Pl}
    |_{\gamma}T_{Pl}|_{\gamma} = \tilde{\nabla}_{Pl}|_{\gamma} R_{Pl}
    |_{\gamma}= 0~.
  \end{equation}
  Define $\tilde{\nabla}_{Pl}(u,v,y^i) :=
  \tilde{\nabla}_{Pl}|_{\gamma}(u,0,0)$. We have $g_{Pl}$ is independent
  of $v,y^i$ so
  \begin{equation}
    \tilde{\nabla}_{Pl} g_{Pl} = \tilde{\nabla}_{Pl} T_{Pl} =
    \tilde{\nabla}_{Pl} R_{Pl} = 0~.
  \end{equation}
  Applying Theorem \ref{1} gives the result.
\end{proof}

\begin{cor} A homogeneous structure $T$ has a well--defined
Penrose limit along a null geodesic $\gamma(t)$ if and only if
$\gamma(t)$ can be re--parameterized to a geodesic of the
canonical connection with respect to $T$.

\end{cor}
\begin{proof} $T$ has a well--defined limit if and only if
$T_{uu}^{y^i}|_{\gamma} = T_{uu}^{v}|_{\gamma}=0$. The proof of
Proposition \ref{red} shows that this is the case if and only if
$\gamma$ can be re--parameterized to a geodesic of the canonical
connection.
\end{proof}

\section{Homogeneous Plane Waves}

We can learn more about the hereditary properties of homogeneity
by studying the space of homogeneous plane waves. In
\cite{Blau:2002js}, Blau and O'Loughlin have classified all
homogeneous plane waves into two classes. The first class consists
of complete metrics and the second class incomplete metrics:

\begin{thm}[Blau--O'Loughlin \cite{Blau:2002js}] \label{homogpw} There are two classes of
  homogeneous plane waves:
  \begin{enumerate}
  \item $g = 2dx^+ dx^- + (e^{x^+f}A_0e^{-x^+f})_{ij}z^iz^j(dx^+)^2 +
    \sum_{i} (dz^i)^2$. Complete metrics.
  \item $g = 2dx^+dx^- +(e^{f \log x^+}A_0e^{-f \log x^+})_{ij}
    z^iz^j \frac{(dx^+)^2}{(x^+)^2} + \sum_{i} (dz^i)^2$. Incomplete
    metrics (singularity along $x^+$).
  \end{enumerate}
\end{thm} The isometry algebra of the generic homogeneous plane wave is
given by:

\begin{eqnarray*} &&[e_i, Y_j] = c \delta_{ij} Z, \qquad [e_i, X]
=-Y_i, \\
&&[Y_i,Y_j] = 2c f_{ij} Z, \qquad [X,Z] = a Z \\
&&[X,Y_i] = (a \delta_{ij} +2 f_{ij} ) Y_j + (c(a + b)^2
(A_0)_{ij}-a f_{ij}- f_{ik}f_{kj}) e_j
\end{eqnarray*} Here $(A_0)_{ij}$ is symmetric and $f_{ij}$ skew--symmetric.
The isotropy is generated by the $e_i$'s. From this it is clear
that homogeneous plane waves are reductive. The non--singular
plane waves have an isometry algebra with $b=c=1$ and $a=0$, while
the singular plane waves have an algebra with $a=c=1$ and $b=0$.
By calculating the homogeneous structure associated to these
reductive splittings we see that the non--singular plane waves are
naturally reductive, while the singular plane waves are not.

Contained in the class of naturally reductive plane waves are the
symmetric plane waves, also
    called the Cahen--Wallach spaces (see \cite{CahenWallach} for the original
    paper or \cite{FOPflux}.) These are given by taking $f_{ij}=0$
    in $(1)$ of Theorem \ref{homogpw}
    and can be diagonalised to the form:
 \begin{equation*}
      g = 2dx^+ dx^- + \sum_{i}A_{i}(z^i)^2(dx^+)^2 +
      \sum_{i} (dz^i)^2
    \end{equation*}
    with $A_{i}$ constant.

Combining this classification with Corollary \ref{plhomogcor} we
obtain the

    \begin{thm} The Penrose limit of a lorentzian space along an
absolutely homogeneous geodesic is a naturally reductive plane
wave.
\end{thm} Also we have the

\begin{prop}
  The Penrose limit of a geodesically complete lorentzian metric $g$
  along a homogeneous geodesic is naturally reductive homogeneous.
\end{prop}

\begin{proof}
  If $g$ is geodesically complete then the Penrose limit is
  complete. The classification of homogeneous plane-waves shows that a
  complete homogeneous plane-wave is naturally reductive.
\end{proof}

\section{Examples}

In this section we will give some examples to show that the above
Theorems cannot be strengthened any further.

First we will show that the converse to Theorem \ref{plhomog} is
not true, i.e. we give an example which shows that the Penrose
limit of a non--homogeneous geodesic in a non--homogeneous space
may be homogeneous. Consider the metric
\begin{equation} g = 2du dv + u dv^2 + \sqrt{u} \sum_i (dx^i)^2.
\end{equation} This is an incomplete and non--homogeneous metric with no
Killing vector in the $\partial_u$ direction. Therefore the null
geodesic given by $\partial_u$ is not homogeneous. However the
Penrose limit of $(g, \partial_u)$ is given by
\begin{equation}  2du dv + \sqrt{u} \sum_i (dx^i)^2.
\end{equation} This is a reductive plane wave \cite{Papadopoulos:2002bg}:

Next we will consider non--homogeneous geodesics in a homogeneous
space. In \cite{Patricot:2003dh} Patricot calculated the Penrose
limits of the \textbf{Kaigorodov space} $K_{n+3}$ which is
$\R^{n+3}$ together with the metric:
\begin{equation*}
  g_{n+3} = e^{-2nL \rho} dx^2 + e^{4L \rho} (2dxdt + \sum_{i=1}^n
  (dy^i)^2) + d \rho^2~,
\end{equation*}
where $L= \frac{1}{2} \sqrt{ - \frac{\Lambda}{n+2}}$.

This is a homogeneous space whose isometries are generated by
\begin{eqnarray*}  &&\boldsymbol{K}_{(0)} = \frac{\partial}{\partial t},
\qquad \boldsymbol{K}_{(x)} = \frac{\partial}{\partial x}, \qquad
\boldsymbol{K}_{(i)} = \frac{\partial}{\partial y^i}, \qquad
\boldsymbol{L}_{i} = x \frac{\partial}{\partial y^i} - y^i
\frac{\partial}{\partial t}, \\  & &\boldsymbol{L}_{ij} = y^i
\frac{\partial}{\partial y^j} - y^j \frac{\partial}{\partial y^i},
\qquad \boldsymbol{J}= \frac{\partial}{\partial \rho} -
at\frac{\partial}{\partial t} - bt \frac{\partial}{\partial x} -
cy^i\frac{\partial}{\partial y^i}, \nonumber\\ \end{eqnarray*}
Here $a = (n+4)L$, $b =-nL$ and $c =2L$.

\begin{footnotesize}
  \begin{center}
    \begin{tabular}{>{$}c<{$}|>{$}c<{$}>{$}c<{$}>{$}c<{$}>{$}c<{$}>{$}c<{$}>{$}c<{$}}
      [,] & \boldsymbol{L}_r & \boldsymbol{L}_{uv} & \boldsymbol{K}_0 & \boldsymbol{K}_x &
      \boldsymbol{K}_i & \boldsymbol{J}\\
      \hline
      \boldsymbol{L}_r & 0 & \delta_{ur} \boldsymbol{K}_v -\delta_{vr}
      \boldsymbol{K}_u &0 & -\boldsymbol{K}_r & \delta_{ri} \boldsymbol{K}_0 & (c-a)
      \boldsymbol{L}_r \\
      \boldsymbol{L}_{st} & -\delta_{sr} \boldsymbol{L}_t + \delta_{tr}
      \boldsymbol{L}_s & \boldsymbol{L}_{stuv} & 0 & 0 & \delta_{ti}
      \boldsymbol{K}_s - \delta_{is} \boldsymbol{K}_t & 0 \\
      \boldsymbol{K}_0 & 0 & 0 & 0 & 0 & 0 & -a\boldsymbol{K}_0 \\
      \boldsymbol{K}_x & \boldsymbol{K}_r & 0 & 0 & 0 & 0 & -b\boldsymbol{K}_x \\
      \boldsymbol{K}_i & -\delta_{ri} \boldsymbol{K}_0 & \delta_{iu}\boldsymbol{K}_v
      - \delta_{iv}\boldsymbol{K}_u & 0 & 0 & 0 & -c\boldsymbol{K}_i \\
      \boldsymbol{J} & (a-c)\boldsymbol{L}_r & 0 & a\boldsymbol{K}_0 & b\boldsymbol{K}_x &
      c\boldsymbol{K}_i & 0
    \end{tabular}
  \end{center}
\end{footnotesize}
where $\boldsymbol{L}_{stuv} = -\delta_{su} \boldsymbol{L}_{tv} +
\delta_{tu} \boldsymbol{L}_{sv} -\delta_{jl} \boldsymbol{L}_{ik} +
\delta_{il} \boldsymbol{L}_{jk} $. This isometry algebra is the
semidirect product of an extended Heisenberg algebra and
$\mathfrak{so}(n)$ and has dimension $ (2n + 3)+
\frac{1}{2}n(n+1)$. The homogeneous space as given by the full
isometry group is non--reductive. However, the transitive subgroup
which is generated by $\boldsymbol{K}_0, \boldsymbol{K}_x,
\boldsymbol{K}_i$ and $\boldsymbol{J}$ gives $K_{n+3}$ as a
reductive homogeneous space.

$K_{n+3}$ has $3$ non-isometric Penrose limits. The first is along
a Killing vector and is thus flat. The second along a homogeneous
geodesic is a reductive homogeneous plane wave and the third is
along a non--homogeneous geodesic is non--homogeneous.

Patricot also considered the space $K_{n+3} \times S^d$ where
$S^d$ is sphere with round metric which is again a non-reductive
homogeneous space. It has four non-isometric Penrose limits. Three
along geodesics which are constant on the sphere and hence have
the same Penrose limits as $K_{n+3}$; the flat metric and a
reductive plane wave and a non--homogeneous metric. The fourth
Penrose limit is along a non-homogeneous geodesic which wraps
around the sphere and $K_{n+3}$ and is also a non-homogeneous
plane wave.

Thus the Penrose limit of a non-reductive homogeneous space along
a non-homogeneous geodesic is not necessarily homogeneous. The
next example will illustrate the existence of non--canonical
homogeneous geodesics and will show that the Penrose limit of a
homogeneous structure along such a curve will blow up.

B.~Komrakov~Jnr has compiled a complete classification of
$4$-dimensional pseudo-riemannian homogeneous spaces
\cite{komrakov}.  In his paper he considers the isotropy
representation $\rho: \mathfrak{h} \to \mathfrak{gl(g/h)}$ of a
homogeneous space $G/H$ and classifies first all the complex forms
and then the real forms of the subalgebra $(\rho( \mathfrak{h})
)^{\C} \subset \mathfrak{ so} (4, \C)$. The result is a list of
the possible Lie algebras $\mathfrak{g}$ and chosen subalgebras
$\mathfrak{h}$ and the associated isotropy representation given as
a matrix.

We can then use the Maurer--Cartan form to recover the metric from
$B$. We summarize below some of the properties of this
classification:

\begin{itemize}

\item Number of isotropy representations admitting Riemannian metrics:
  $6$

\item Number of isotropy representations admitting lorentzian
  metrics:$14$

\item Number of isotropy representations admitting metrics of $(2,2)$
  signature: $30$

(There is some overlap in these cases where a representation
admits metrics of different signatures.)

\item Number of symmetric/reductive algebras admitting a Riemannian
  metric: $21$/$29$

\item Number of symmetric/reductive/nonreductive algebras admitting a
 lorentzian metric: $35$/$64$/$6$

\item Number of symmetric/reductive/non-reductive algebras admitting
  a metric of $(2,2)$ signature: $57$/$123$/$9$
\end{itemize}

By studying Komrakov's list we see that there does not exist a
$4$--dimensional Lorentzian homogeneous space with a
non--canonical homogeneous geodesic. However there do exist
$5$--dimensional examples as we will now show. Consider the
algebra (Komrakov number $1.1^2.11$ extended by a central
element.)

\begin{equation*}
  \begin{tabular}{>{$}c<{$}|>{$}c<{$}>{$}c<{$}>{$}c<{$}>{$}c<{$}>{$}c<{$}>{$}c<{$}}
    [,] & \boldsymbol{e}_1 & \boldsymbol{u}_1 & \boldsymbol{u}_2 &
    \boldsymbol{u}_3 & \boldsymbol{u}_4 & \boldsymbol{u}_5\\
    \hline
    \boldsymbol{e}_1 & 0 & \boldsymbol{u}_3&\frac{1}{2}\boldsymbol{u}_4 & -\boldsymbol{u}_1 & -\frac{1}{2}\boldsymbol{u}_2 &0\\
    \boldsymbol{u}_1 & -\boldsymbol{u}_3 & 0 & \boldsymbol{u}_2 & -4 \boldsymbol{e}_1 & -\boldsymbol{u}_4 &0 \\
    \boldsymbol{u}_2 & -\frac{1}{2} \boldsymbol{u}_4 & \boldsymbol{u}_2 & 0 & -\boldsymbol{u}_4 & 0 &0\\
    \boldsymbol{u}_3 & \boldsymbol{u}_1 & 4 \boldsymbol{e}_1 & \boldsymbol{u}_4 & 0 &
    \boldsymbol{u}_2&0\\
    \boldsymbol{u}_4 & \frac{1}{2} \boldsymbol{u}_2 & \boldsymbol{u}_4 & 0 &
    -\boldsymbol{u}_2 & 0&0\\
    \boldsymbol{u}_5 &0&0&0&0&0&0
  \end{tabular}\vspace{1em}
\end{equation*} This defines a reductive homogeneous space $G/H$ with $\mathfrak{m}$ the span of $\{\boldsymbol{u}_1,
\boldsymbol{u}_2,
    \boldsymbol{u}_3 , \boldsymbol{u}_4 , \boldsymbol{u}_5 \}$ and $\mathfrak{h}$ spanned by $\boldsymbol{e}_1$.
    The corresponding isotropy representation is skew--symmetric with respect to the bilinear
form $B$:
\begin{equation}
  \begin{pmatrix}
    1 & 0 & 0 & 0 &0 \\
    0 & 1 & 0 & 0 &0\\
    0 & 0 & 1 & 0 &0\\
    0 & 0 & 0 & 1 &0\\
    0 & 0 & 0 & 0 & -1
  \end{pmatrix}
\end{equation} To determine the induced metric we make a choice of
local coset representative
\begin{equation} \sigma = \exp(x_1\boldsymbol{u}_1) \exp(x_2
  \boldsymbol{u}_2) \exp(x_3 \boldsymbol{u}_3) \exp(x_4 \boldsymbol{u}_4) \exp(x_5 \boldsymbol{u}_5):
  M \to G~,
  \end{equation} and calculate the Maurer--Cartan form

  \begin{eqnarray*}  \sigma^{-1} d \sigma &= &
  \cosh(2x_3)dx_1\boldsymbol{u}_1 \\ &&+ (x_4\sinh(2x_3)dx_1 +x_2\cosh(x_3)dx_1
  +\cosh(x_3)dx_2+x_4dx_3)\boldsymbol{u}_2 \\ &&+ dx_3 \boldsymbol{u}_3\\
  &&
   + (-x_4 \cosh(2x_3)dx_1-2
  \sinh(2x_3)dx_1-x_3 \sinh(x_3)dx_1-\sinh(x_3)dx_2+dx_4)
  \boldsymbol{u}_4 \\ &&+ dx_5 \boldsymbol{u}_5. \end{eqnarray*}
  The metric is given by $B((\sigma^{-1} d \sigma)_{\mathfrak{m}}, (\sigma^{-1} d
  \sigma)_{\mathfrak{m}})$. The non--zero components of the homogeneous
  structure $T$ restricted to the subspace $\mathfrak{m}$ are given by:

  \begin{eqnarray*}
&&1=T_{414}= T_{221}=T_{243}=T_{244}=T_{423}=T_{424}, \\
&&-1=T_{212}=T_{234}= T_{432}= T_{441}.
\end{eqnarray*}

Now consider the vector $U = \boldsymbol{u}_2 + \frac{1}{\sqrt{2}}
\boldsymbol{u}_3+ \sqrt{\frac{3}{2}}\boldsymbol{u}_5 + \sqrt{2}
\boldsymbol{e}_1$. This is an absolutely geodetic vector and hence
generates a homogeneous geodesic. However this geodesic is not a
geodesic of the canonical connection and thus the Penrose limit of
the homogeneous structure along $\exp( tU)(p)$ will blow up:
\begin{equation} T(U,U)|_{\mathfrak{m}} = \boldsymbol{u}_1 - \frac{1}{2}
\boldsymbol{u}_4 \xrightarrow{Pl} \infty. \end{equation}

However, the algebra $\mathfrak{g}$ has the following transitive
subalgebra:

\begin{equation*}
  \begin{tabular}{>{$}c<{$}|>{$}c<{$}>{$}c<{$}>{$}c<{$}>{$}c<{$}>{$}c<{$}}
    [,] & U & \boldsymbol{u}_1 & \boldsymbol{u}_2 &
    \boldsymbol{u}_4 & \boldsymbol{u}_5\\
    \hline
    U & 0 & 2U- 3\boldsymbol{u}_2&\sqrt{2} \boldsymbol{u}_4 & 0 &0\\
    \boldsymbol{u}_1 & -2U+ 3\boldsymbol{u}_2 & 0 & \boldsymbol{u}_2 & - \boldsymbol{u}_4 &0 \\
    \boldsymbol{u}_2 & \sqrt{2} \boldsymbol{u}_4 & -\boldsymbol{u}_2 & 0 & 0 &0\\
    \boldsymbol{u}_4 & 0 &  \boldsymbol{u}_4 & 0 & 0 &0\\
    \boldsymbol{u}_5 &0&0&0&0&0
  \end{tabular}\vspace{1em}
\end{equation*} The non--zero components of the homogeneous structure are
given by
 \begin{eqnarray*}
&&-1 = T_{U12}=T_{441},\\ &&1=T_{221}=T_{21U}=T_{212}=T_{414},\\
&& \frac{1}{\sqrt{2}}=T_{U24}=T_{24U}.
\end{eqnarray*} This homogeneous structure
does not blow up under the Penrose limit along $U$. It is not
clear that every homogeneous geodesic is canonical with respect to
some reductive decomposition as in this case. However, this is
true for all the null homogeneous geodesics of $4$--dimensional
lorentzian homogeneous spaces.

Finally we will show that we can not replace homogenous with
reductive in the result on Penrose limits of lorentzian metrics
along homogeneous geodesics. Consider the incomplete, nonreductive
plane wave metric
\begin{equation*}
  dudv + u^{2 \mu} (dy^i)^2.
\end{equation*}
This is part of Blau and O'Loughlin's classification of
homogeneous plane waves. It has a singularity at $u=0$ for $\mu
\neq 0$. The vector field
\begin{equation*}
  X = -u \frac{\partial}{\partial u} + v
  \frac{\partial}{\partial v}   + 2 \mu y^i \frac{\partial}{\partial
    y^i}
\end{equation*}
is a Killing vector and thus the geodesic defined by $\partial_u$
is a homogeneous geodesic. The trivial Penrose limit along this
geodesic gives the same metric and therefore shows that the
Penrose limit of a non-reductive space along a homogeneous
geodesic is not necessarily reductive.

\section{The existence of homogeneous geodesics}

Finally we would like to make a remark on the existence of null
homogeneous geodesics. The following Theorem has been proven in
\cite{MR2002b:53076}, \cite{MR2001f:53104}.

\begin{thm}[Kowalski--Szenthe]
  Every homogeneous Riemannian manifold admits at least one
  homogeneous geodesic through every point.
\end{thm}

Since every homogenous Riemannian manifold is reductive (see
\cite{MR85b:53052}) it appears that this Theorem is also true in
the case of reductive lorentzian manifolds.

\begin{prop}
  Every reductive homogeneous lorentzian manifold admits at least one
  homogeneous geodesic through every point.
\end{prop} In fact all lorentzian homogeneous examples known to the author (and this includes all
$4$--dimensional homogeneous spaces appearing on Komrakov's list,)
contain at least one null homogeneous geodesic although not all of
them contain an absolutely homogeneous one as the following
example (Komrakov number $1.1^2$) shows:
\begin{table*}[h!]
  \centering
  \begin{tabular}{>{$}c<{$}|>{$}c<{$}>{$}c<{$}>{$}c<{$}>{$}c<{$}>{$}c<{$}}
    [,] & \boldsymbol{e}_1 & \boldsymbol{u}_1 & \boldsymbol{u}_2 & \boldsymbol{u}_3 &
    \boldsymbol{u}_4 \\
    \hline
    \boldsymbol{e}_1 & 0 & \boldsymbol{u}_3&0 & -\boldsymbol{u}_1 & 0 \\
    \boldsymbol{u}_1 & -\boldsymbol{u}_3 & 0 & 0 & -\boldsymbol{u}_2 & \boldsymbol{u}_1 \\
    \boldsymbol{u}_2 & 0 & 0 & 0 & 0 & 2\boldsymbol{u}_2 \\
    \boldsymbol{u}_3 & \boldsymbol{u}_1 & \boldsymbol{u}_2 & 0 & 0 &  \boldsymbol{u}_3 \\
    \boldsymbol{u}_4 & 0 & -\boldsymbol{u}_1 & -2\boldsymbol{u}_2 & -\boldsymbol{u}_3 & 0
  \end{tabular}
  \label{tab:1.1-2}
\end{table*}
together with the bilinear form
\begin{equation*}
  B =
  \begin{pmatrix}
    1 & 0 & 0 & 0 \\
    0 & 1 & 0 & 0 \\
    0 & 0 & 1 & 0 \\
    0 & 0 & 0 & -1
  \end{pmatrix}
\end{equation*}
This is a reductive algebra and so using Proposition \eqref{geo}
it can be shown that the homogeneous space derived from this
algebra and bilinear form has no null absolutely homogeneous
geodesics. (This is in fact effectively the only $4$-dimensional
lorentzian homogeneous space without any null absolutely
homogeneous geodesics.) However it does have a family of null
homogeneous geodesics.

\begin{eqnarray*} &&U = A\boldsymbol{u}_4 \pm A\boldsymbol{u}_2 + B
\boldsymbol{e}_1 \textrm{ and } \lambda = -2A \quad \textrm{ with } A,B \in \R \\
&&\textrm{or}\\
&&U = A \boldsymbol{u}_2 + B \boldsymbol{u}_3 + C\boldsymbol{u}_4
\textrm{ and } \lambda = -C \quad \textrm{ with } A^2 + B^2 = C^2,
\quad A,B,C \in \R.
\end{eqnarray*}

To the author's knowledge there are no known results about the
existence of homogeneous geodesics in the nonreductive case.

\begin{section}{Acknowledgments}
We are very grateful to Jose Figueroa-O'Farrill for help and
encouragement throughout the research and writing up of this
paper. We would also like to thank Patrick Meesseen for his
invaluable correspondence, Harry Braden for drawing our attention
to \cite{MR85b:53052} and to Jos\'e Antonio Oubi\~na for early
useful correspondence.

This research was funded by an EPSRC Postgraduate studentship.

\end{section}

\bibliography{pen}
\bibliographystyle{plain}

\end{document}